\documentclass[12pt,reqno]{amsart}

    \usepackage{amssymb}
    \usepackage{graphicx}
    \usepackage{pstricks,pst-plot}

    \newtheorem{thm}{Theorem}[section]
    \newtheorem{prop}[thm]{Proposition}
    \newtheorem{lem}[thm]{Lemma}

    \theoremstyle{definition}
        \newtheorem{defn}[thm]{Definition}
        \newtheorem{exgr}[thm]{Example}

        \newtheorem{alg}[thm]{Algorithm}
        \newtheorem{game}[thm]{Game}

    \theoremstyle{remark}
        \newtheorem{rem}[thm]{Remark}

    \numberwithin{equation}{section}
    \numberwithin{figure}{section}

    \title[Combinatorics of Patience Sorting Piles]{Combinatorics of Patience Sorting Piles*}

    \author[Alexander Burstein and Isaiah Lankham]{Alexander
    Burstein$^\dagger$ and Isaiah Lankham$^\ddagger$}

    \thanks{* This material is based in part upon work supported by the
    National Science Foundation under Grant No. DMS-0502858.}

    \thanks{$^\dagger$
    Department of Mathematics,
    Iowa State University,
    Ames, IA 50011-2064, USA.
    (\texttt{burstein@math.iastate.edu})
    Research supported in part by the
    U.S. NSA Young Investigator Grant H98230-06-1-0037.
    }

    \thanks{$^\ddagger$
    Department of Mathematics,
    University of California-Davis,
    Davis, CA 95616-8633, USA.
    (\texttt{issy@math.ucdavis.edu})
    Research supported in part by the
    U.S. NSF Grants DMS-0135345 and DMS-0304414.
    }

    \keywords{Patience Sorting, Set Partitions, Bell Numbers,
    Generalized Permutation Patterns, Left-to-right Minima
    Subsequences, Basic Subsequences, Shadow Diagrams}

\begin{document}

    \begin{abstract}
        Despite having been introduced in 1962 by C.L. Mallows, the
        combinatorial algorithm \emph{Patience Sorting} is only now
        beginning to receive significant attention due to such recent
        deep results as the Baik-Deift-Johansson Theorem that connect
        it to fields including Probabilistic Combinatorics and Random
        Matrix Theory.

        The aim of this work is to develop some of the more basic
        combinatorics of the Patience Sorting Algorithm.  In
        particular, we exploit the similarities between Patience
        Sorting and the Schensted Insertion Algorithm in order to do
        things that include defining an analog of the Knuth relations
        and extending Patience Sorting to a bijection between
        permutations and certain pairs of set partitions.  As an
        application of these constructions we characterize and
        enumerate the set $S_{n}(3\textrm{-}\bar{1}\textrm{-}42)$ of
        permutations that avoid the generalized permutation pattern
        $2\textrm{-}31$ unless it is part of the generalized pattern
        $3\textrm{-}1\textrm{-}42$.\\

        \noindent {\sc R\'esum\'e}.  En d\'{e}pit de la introduction
        en 1962 par C.L. Mallows, combinatoire d'al\-gor\-ith\-me
        \emph{Patience Sorting} commence seulement maintenant \`{a}
        susciter l'attention significative d\^{u} \`{a} des
        r\'{e}sultats profonds r\'{e}cents tels que le
        th\'{e}or\`{e}me de Baik-Deift-Johansson qui le relient \`{a}
        la combinatoire probabiliste et \`{a} la th\'{e}orie des
        matrices al\'{e}atoires.

        On d\'{e}veloppe une partie plus fondamentale de la
        combinatoire de l'algorithme de \emph{Patience Sorting}.  En
        particulier, on utilise les similitudes entre \emph{Patience
        Sorting} et la correspondence de Schensted pour d\'{e}finir un
        analogue des relations de Knuth et pour g\'{e}n\'{e}raliser
        \emph{Patience Sorting} en une bijection entre les
        permutations et certaines paires de partitions d'ensemble.
        Comme application de ces constructions on caract\'{e}rise et
        \'{e}num\'{e}re l'ensemble
        $S_{n}(3\textrm{-}\bar{1}\textrm{-}42)$ de permutations qui
        \'{e}vitent le motif g\'{e}n\'{e}ralis\'{e} $2\textrm{-}31$ de
        permutation \`{a} moins qu'il soit partie du motif
        g\'{e}n\'{e}ralis\'{e} $3\textrm{-}1\textrm{-}42$.

    \end{abstract}

    \maketitle

%     %First page headline in AmS-LaTeX for S\'eminaire Lotharingien de Combinatoire
%     %--first part
%     \thispagestyle{myheadings}
%     \font\rms=cmr8
%     \font\its=cmti8
%     \font\bfs=cmbx8
%     \markright{\its S\'eminaire Lotharingien de
%     Combinatoire \bfs 54A \rms (2005), Article~B54A?\hfill}
%     \def\thepage{}
%     %
%     %

    \section{Introduction}
    \label{sec:Introduction}

    The term \emph{Patience Sorting} was introduced in 1962 by
    C.~L.~Mallows \cite{refMallows1962, refMallows1963} as the name of
    a card sorting algorithm invented by A.~S.~C.~Ross.  This
    algorithm works by first partitioning a shuffled deck of cards
    (which throughout this paper we take to be a permutation $\sigma
    \in \mathfrak{S}_{n}$) into sorted subsequences called
    \emph{piles} using what Mallows referred to as a ``patience
    sorting procedure'':

    \begin{alg}[Mallows' Patience Sorting Procedure]
    \label{alg:MallowsPSprocedure} Given a shuffled deck of
    cards $\sigma = c_{1} c_{2} \cdots c_{n}$, inductively build the
    set of piles $R = R(\sigma) = \{r_{1}, r_{2}, \ldots, r_{m}\}$ as
    follows:\smallskip

        \begin{itemize}
            \item Place the first card $c_{1}$ from the deck into a
            pile $r_{1}$ by itself.\smallskip

            \item For each remaining card $c_{i}$ ($i = 2, \ldots,
            n$), consider the cards $d_{1}, d_{2}, \ldots, d_{k}$ atop
            the piles $r_{1}, r_{2}, \ldots, r_{k}$ that have already
            been formed.\smallskip

                \begin{itemize}
                    \item If $c_{i} > \max\{d_{1}, d_{2}, \ldots,
                    d_{k}\}$, then put $c_{i}$ into a new right-most
                    pile $r_{k+1}$ by itself.\smallskip

                    \item Otherwise, find the left-most card $d_{j}$
                    that is larger than $c_{i}$ and put the card
                    $c_{i}$ atop pile $r_{j}$.
                \end{itemize}
        \end{itemize}
    \end{alg}

%     %First page headline in AmS-LaTeX for S\'eminaire Lotharingien de Combinatoire
%     %--restoring the headers and pagenumbering
%     \pagenumbering{arabic}
%     \addtocounter{page}{1}
%     \markboth{\SMALL ALEXANDER BURSTEIN AND ISAIAH LANKHAM}{\SMALL COMBINATORICS OF PATIENCE SORTING PILES}
%     %
%     %

    \noindent We call the collection of piles $R(\sigma)$ the
    \emph{pile configuration} associated to the deck of cards $\sigma
    \in \mathfrak{S}_{n}$ and illustrate their formation via an
    extended version of Algorithm \ref{alg:MallowsPSprocedure} in
    Section \ref{sec:ExtendingPS:Algorithm} below.\smallskip

    Since each card $c_{i}$ is either larger than the top card of
    every pile or is placed on top of the left-most top card $d_{j}$
    larger than it, the cards $d_{1}, d_{2}, \ldots, d_{k}$ atop the
    piles will be in increasing order from left to right at each step
    of the algorithm.  Thus, Algorithm~\ref{alg:MallowsPSprocedure}
    resembles repeated application of the Schensted Insertion
    Algorithm (as discussed in \cite{refAD1999}) for interposing a
    value into the increasing sequence $d_{1}, d_{2}, \ldots, d_{k}$
    as if it were the top row of a Young tableau.  The distinction is
    that cards remain in place and have other cards placed on top of
    them instead of being actively ``bumped'' from the row so that the
    Schensted Insertion Algorithm can then be recursively applied to
    the ``bumped'' value and the next lower row in the Young tableau.
    In this sense, Patience Sorting can be viewed as a non-recursive
    analog of the remarkable Robinson-Schensted-Knuth (or RSK)
    Algorithm due to G.~Robinson \cite{refRobinson1938} for
    permutations in 1938, C.~Schensted \cite{refSchensted1961} for
    words in 1961, and D.~Knuth \cite{refKnuth1970} for so-called
    $\mathbb{N}$-matrices in 1970.  (See Fulton \cite{refFulton1997}
    for the appropriate definitions and for a detailed account of the
    differences between these algorithms.)

    Recall that the RSK Algorithm bijectively associates an ordered
    pair of standard Young tableaux $(P(\sigma), Q(\sigma))$ to each
    permutation $\sigma = \sigma_{1}\sigma_{2} \cdots \sigma_{n} \in
    \mathfrak{S}_{n}$ by first building a so-called ``insertion
    tableau'' $P(\sigma)$ through repeated Schensted Insertion of the
    components $\sigma_{1}, \sigma_{2}, \ldots, \sigma_{n}$ into an
    initially empty tableau.  It also simultaneously constructs the
    ``recording tableau'' $Q(\sigma)$ by literally \emph{recording}
    how $P(\sigma)$ is formed.  These tableaux have the same shape (a
    partition $\lambda$ of $n$, denoted $\lambda \vdash n$), and this
    correspondence has many interesting properties.  E.g., RSK applied
    to a permutation is symmetric in the sense that if $\sigma \in
    \mathfrak{S}_{n}$ corresponds to the ordered pair of tableaux
    $(P(\sigma), Q(\sigma))$, then $(Q(\sigma), P(\sigma))$
    corresponds to the inverse permutation $\sigma^{-1}$.  As a
    result, there is a bijection between the set of involutions
    $\mathfrak{I}_{n} \subset \mathfrak{S}_{n}$ and the set
    $\mathfrak{T}_{n}$ of all standard Young tableaux with entries $1,
    2, \ldots, n$.  (This is the famous Sch\"{u}tzenberger symmetry
    property first proven in \cite{refSchutzenberger1963}.)

    In this paper we develop a bijective extension of Algorithm
    \ref{alg:MallowsPSprocedure} and then study analogues for such
    properties of RSK. To facilitate this, we first characterize in
    Section~\ref{sec:PileConfigurations} when two permutations have
    the same pile configurations under Algorithm
    \ref{alg:MallowsPSprocedure}.  This yields an equivalence relation
    $\stackrel{PS}{\sim}$ on $\mathfrak{S}_{n}$ that is analogous to
    the Knuth relation $213 \stackrel{RSK}{\sim} 231$.  (Recall that
    the Knuth relations describe when two permutations have the same
    ``insertion tableau'' $P$ under RSK; see Sagan
    \cite{refSagan2000}.)

    In Section \ref{sec:ExtendingPS} we then explicitly describe a
    bijection between $\mathfrak{S}_{n}$ and certain pairs of pile
    configurations having the same shape (a composition $\gamma$ of
    $n$, denoted $\gamma \models n$).  Since there are many more
    possible pile configurations than standard Young tableaux (the
    former are enumerated by Bell numbers; see Theorem
    \ref{thm:EnumeratingSn3142}), it is necessary to specify which
    pairs are possible; this turns out to involve the same patterns as
    the other Knuth relation $312 \stackrel{RSK}{\sim} 132$ (see
    Definition \ref{defn:StablePairs}).  Moreover, this bijection
    shares the same Sch\"{u}tzenberger symmetry property as RSK, and
    so we can immediately characterize a certain collection of pile
    configurations that are in bijection with the set of involutions
    $\mathfrak{I}_{n}$ (as well as with the set $\mathfrak{T}_{n}$ of
    standard Young tableaux).

    In Section \ref{sec:EnumeratingSn3142} we conclude by using the
    equivalence relation $\stackrel{PS}{\sim}$ to characterize and
    enumerate the set $S_{n}(3\textrm{-}\bar{1}\textrm{-}42)$ of
    permutations avoiding the barred (generalized) permutation pattern
    $3\textrm{-}\bar{1}\textrm{-}42$.  Such permutations avoid the
    pattern $2\textrm{-}31$ unless it is contained in a
    $3\textrm{-}1\textrm{-}42$ pattern.  (See Sections
    \ref{sec:PileConfigurations:PSEquivalence} and
    \ref{sec:EnumeratingSn3142} for the appropriate
    definitions.)\bigskip

    Another interesting property of RSK is that, given $\sigma \in
    \mathfrak{S}_{n}$, the number of boxes in the top row of the
    ``insertion tableau'' $P(\sigma)$ is exactly the length of the
    longest increasing subsequence in $\sigma$.  (This was first
    proven by Schensted \cite{refSchensted1961} but is now a special
    case of Greene's Theorem \cite{refGreene1974}).  Due to the
    similarity between the Schensted Insertion Algorithm and Algorithm
    \ref{alg:MallowsPSprocedure}, it is clear that the cards atop the
    piles when Patience Sorting terminates will be exactly the
    elements in the top row of $P(\sigma)$.  Thus, the number of piles
    formed under Patience Sorting is also equal to the length of the
    longest increasing subsequence in $\sigma$, and so one can apply
    the recent but now highly celebrated Baik-Deift-Johansson Theorem
    \cite{refBDJ1999} in order to describe the asymptotic distribution
    for the number of piles (up to rescaling).  Due to this deep
    connection between Patience Sorting and Probabilistic
    Combinatorics, it has been suggested (see, e.g.,
    \cite{refKuy1999}, \cite{refMackenzie1998} and
    \cite{refPeterson1999}; cf.  \cite{refConrey2003}) that studying
    generalizations of Patience Sorting might be the key to tackling
    certain open problems that can be viewed from the standpoint of
    Random Matrix Theory---the most notable being the Riemann
    Hypothesis.

    At the same time, there is a lot more to Patience Sorting than
    just resembling the RSK Algorithm for permutations.  E.g., after
    applying Algorithm \ref{alg:MallowsPSprocedure} to a deck of
    cards, it is easy to recollect each card in ascending order from
    amongst the current top cards of the piles (and thus complete
    A.~S.~C.~Ross' card sorting algorithm).  While this is not
    necessarily the fastest sorting algorithm one can apply to a deck
    of cards, the \emph{patience} in \emph{Patience Sorting} is not
    intended to describe a prerequisite for its use.  Instead it
    refers to how pile formation in Algorithm
    \ref{alg:MallowsPSprocedure} resembles the way in which one places
    cards into piles when playing the popular single-person card game
    \emph{Klondike Solitaire}, which is often called \emph{Patience}
    in the UK. This is more than a coincidence, though, as Algorithm
    \ref{alg:MallowsPSprocedure} also happens to be an optimal
    strategy (in the sense of forming as few piles as possible; see
    \cite{refAD1999} for a proof) when playing an idealized model of
    Klondike Solitaire known as \emph{Floyd's Game}:

    \begin{game}[Floyd's Game]\label{game:FloydsGame} Given a shuffled
    deck of cards $c_{1}, c_{2}, \ldots, c_{n}$,

        \begin{itemize}
            \item Place the first card $c_{1}$ from the deck into a
            pile $r_{1}$ by itself.

            \item Then, for each card $c_{i}$ ($i = 2, \ldots, n$),
            either

                \begin{itemize}

                    \item put $c_{i}$ into a new pile by itself or

                    \item play $c_{i}$ on top of any pile whose
                    current top card is larger than $c_{i}$.

                \end{itemize}

            \item The object of the game is to end with as few piles
            as possible.

        \end{itemize}

    \end{game}

    \noindent In other words, the cards are played one at a time
    according to the order they appear in the deck so that piles are
    created in much the same way they are formed under Patience
    Sorting.  According to \cite{refAD1999}, Floyd's Game was
    developed independently of Mallow's work and originated in
    unpublished correspondence between computer scientists Bob Floyd
    and Donald Knuth during 1964.

    Note that unlike Klondike Solitaire, there is a known strategy
    (Algorithm \ref{alg:MallowsPSprocedure}) for Floyd's Game under
    which one will always win.  In fact, Klondike Solitaire---though
    so popular that it has come pre-installed on the vast majority of
    personal computers shipped since 1989---is very poorly understood
    mathematically.  (Recent progress, however, has been made in
    developing an optimal strategy for a version called
    \emph{thoughtful solitaire} \cite{refYDRR2005}.)  As such, Persi
    Diaconis (\cite{refAD1999} and private communication with the
    second author) has suggested that a deeper understanding of
    Patience Sorting and its generalization would undoubtedly help in
    developing a better mathematical model for analyzing Klondike
    Solitaire.

    \section{Pile Configurations Coming from Patience Sorting}
    \label{sec:PileConfigurations}

        \subsection{Pile Configurations, Shadow Diagrams, and Reverse Patience Words}
        \label{sec:PileConfigurations:Characterization}

        We begin by explicitly characterizing the pile configurations
        that result from applying Patience Sorting (Algorithm
        \ref{alg:MallowsPSprocedure}) to a permutation:

        \begin{lem}\label{lemma:PileConfiguration}
            Let $\sigma \in \mathfrak{S}_{n}$ be a permutation and
            $R(\sigma) = \{r_{1}, r_{2}, \ldots, r_{k}\}$ be the pile
            configuration associated to $\sigma$ under Algorithm
            \ref{alg:MallowsPSprocedure}.  Then $R(\sigma)$ is a
            partition of the set $[n] = \{1, 2, \dots, n\}$ such that
            denoting $r_{j} = \{ r_{j 1}> r_{j 2} > \dots > r_{j
            s_{j}}\}$,
                \begin{equation}\label{eq:PileConfigurationCondition}
                    r_{j s_{j}} < r_{i s_{i}} \quad \mathrm{if} \quad j < i
                \end{equation}

            \noindent Moreover, for every set partition $T = \{t_{1},
            t_{2}, \ldots, t_{k}\}$ satisfying Equation
            (\ref{eq:PileConfigurationCondition}), there is a
            permutation $\tau \in \mathfrak{S}_{n}$ such that $R(\tau)
            = T$.

        \end{lem}

        \begin{proof}
        Given a pile configuration $R(\sigma) = \{r_{1}, r_{2},
        \ldots, r_{k}\}$, suppose that for some pair of indices
        $i,j\in [k]$ we have $j<i$ but $r_{j s_{j}} > r_{i s_{i}}$.
        Then $r_{i s_{i}}$ was put atop pile $r_i$ when pile $r_j$ had
        top card $d_{j} \geq r_{j s_{j}}$ so that $d_{j} > r_{i
        s_{i}}$.  However, it then follows that, under
        Algorithm~\ref{alg:MallowsPSprocedure}, the card $r_{i s_{i}}$
        would actually then have been placed atop either pile $r_j$ or
        a pile to the left of $r_j$ instead of atop pile $r_i$.  The
        resulting contradiction implies that $r_{j s_{j}} < r_{i
        s_{i}}$ for each $j<i$.

        Conversely, let $T = \{t_{1}, t_{2}, \dots, t_{k}\}$ be any
        set partition of $[n]$ with $t_{j} = \{ t_{j 1}> t_{j 2} >
        \dots > t_{j s_{j}} \}$ for every $j\in[k]$ and $t_{j s_{j}} <
        t_{i s_{i}}$ for all pairs $j < i$.  Consider then the
        permutation
        \[
        \tau= t_{1 1} t_{1 2} \dots t_{1 s_{1}} \, t_{2 1} t_{2 2}
        \dots t_{2 s_{2}} \; \dots \; t_{k 1} t_{k 2} \dots t_{k
        s_{k}} \in \mathfrak{S}_{n}.
        \]
        We show that $R(\tau)=T$: Given any $j\in[k]$ and any $i>j$,
        each $m\in[s_i]$ satisfies $t_{j s_j} < t_{i s_i} < t_{im}$,
        and so no value to the right of $t_{j s_j}$ in $\tau$ will be
        placed on top of $t_{j s_j}$ in $R(\tau)$.  Since $t_{1
        1}<t_{1 2}<\dots<t_{1 s_{1}}$, all these entries must be
        placed in the leftmost (a.k.a. first) pile $r_{1} \in
        R(\tau)$, so the first pile of $R(\tau)$ is $r_{1} = t_1$.
        Now suppose that the entries $t_{1}$ through $t_{j s_j}$ of
        $\tau$ have been placed in piles $t_1,t_2,\dots, t_j$.  Then
        the entries $t_{j+1,1}>t_{j+1,2}>\dots>t_{j+1,s_{j+1}}$ cannot
        be placed atop any $t_{i s_i}$ for $i\le j$, so we must form
        at least one new pile, say $t_{j+1}$.  However, since these
        entries occur in decreasing order in $\tau$, these cards will
        then all be placed in the $(j+1)^{\textrm{st}}$ pile
        $t_{j+1}$.  Moreover, no entry to the right of
        $t_{j+1,s_{j+1}}$ can be placed in pile $t_{j+1}$, so the
        $(j+1)^{\textrm{st}}$ pile of $R(\tau)$ is $r_{j+1} =
        t_{j+1}$.  Therefore, we obtain $R(\sigma)=T$ by induction.
        \end{proof}

        We will often express a pile configuration $R$ with its
        constituent piles $r_{1}, r_{2}, \ldots, r_{k}$ written
        vertically and bottom-justified with respect to the largest
        value $r_{j 1}$ in each pile $r_{j}$.  This motivates the
        following definition (which reverses the so-called
        ``far-eastern reading''):

        \begin{defn}\label{defn:ReversePatienceWords}
            The \emph{reverse patience word} $RPW(R)$ for a pile
            configuration $R$ is the permutation formed by
            concatenating the piles $r_{1}, r_{2}, \ldots, r_{k}$
            together with each pile $r_{j}$ written in decreasing
            order (i.e., read from bottom to top in order from left to
            right).  In the notation of Lemma
            \ref{lemma:PileConfiguration},
            \[
            RPW(R) = r_{1 1} r_{1 2} \dots r_{1 s_{1}} \
            r_{2 1} r_{2 2} \dots r_{2 s_{2}} \quad \dots \quad
            r_{k 1} r_{k 2} \dots r_{k s_{k}}.
            \]

        \end{defn}

        \begin{exgr}\label{eg:PileConfigurationExample}
            The pile configuration $R = \{\{6 > 4 > 1\}, \{5 > 2\},
            \{8 > 7 > 3\}\}$ is represented by the piles

            \begin{center}
                \begin{tabular}{l l l}
                    1 &   & 3 \\
                    4 & 2 & 7 \\
                    6 & 5 & 8 \end{tabular}\\
            \end{center}

            \noindent and has the reverse patience word $RPW(R) =
            64152873$.  Moreover, note that as in the proof of Lemma
            \ref{lemma:PileConfiguration}, $R(RPW(R))=R(64152873) =
            R$.
        \end{exgr}

        The following Lemma should now be clear from the above
        definitions and example:

        \begin{lem}
            Given a permutation $\sigma \in \mathfrak{S}_{n}$,
            $R(RPW(R(\sigma)))=R(\sigma)$.
        \end{lem}

        \begin{proof}
            From the proof of Lemma \ref{lemma:PileConfiguration}, we
            have that if $\tau=RPW(T)$, where $T$ is any set partition
            of $[n]$ satisfying condition
            (\ref{eq:PileConfigurationCondition}), then $T=R(\tau)$.
            Thus, $R(RPW(T))=T$ for any partition $T$ of $[n]$
            satisfying (\ref{eq:PileConfigurationCondition}), in
            particular for any $T=R(\sigma)$ where $\sigma\in
            \mathfrak{S}_{n}$.
        \end{proof}

        \begin{figure}[t]\label{fig:ShadowExample}
              \centering

              \begin{tabular}{cccc}

                  \begin{minipage}[c]{1.325in}
                      \centering

                          \psset{xunit=0.125in,yunit=0.125in}

                          \begin{pspicture}(0,0)(9,9)

                              \psaxes[Dy=2]{->}(9,9)

                              \psframe[fillcolor=gray,fillstyle=solid,linecolor=gray](2,4)(9,9)

                              \rput(5.5,6.5){\psframebox*[framearc=.5]{{\large
                              $S(2,4)$}}}

                              \rput(2,4){{\Large $\bullet$}}

                          \end{pspicture}\\[0.25in]

                        (a) The Shadow $S_{NE}(2,4)$.
                  \end{minipage}
                  &
                  \begin{minipage}[c]{1.5in}
                      \centering

                          \psset{xunit=0.125in,yunit=0.125in}

                          \begin{pspicture}(0,0)(9,9)

                              \psaxes[Dy=2]{->}(9,9)

                              \pspolygon[fillcolor=gray, linecolor=gray,
                              fillstyle=solid](1,9)(1,6)(2,6)(2,4)(4,4)(4,1)(9,1)(9,9)

                              \psline[linecolor=darkgray,linewidth=2pt](1,9)(1,6)(2,6)(2,4)(4,4)(4,1)(9,1)

                               \rput(1,6){{\Large $\bullet$}}%
                               \rput(2,4){{\Large $\bullet$}}%
                               \rput(3,5){{\Large $\bullet$}}%
                               \rput(4,1){{\Large $\bullet$}}%
                               \rput(5,8){{\Large $\bullet$}}%
                               \rput(6,7){{\Large $\bullet$}}%
                               \rput(7,2){{\Large $\bullet$}}%
                               \rput(8,3){{\Large $\bullet$}}%

                          \end{pspicture}\\[0.25in]

                      (b) Shadowline $L_{1}(64518723)$.
                  \end{minipage}

                  \begin{minipage}[c]{1.5in}
                      \centering

                          \psset{xunit=0.125in,yunit=0.125in}

                          \begin{pspicture}(0,0)(9,9)

                              \psaxes[Dy=2]{->}(9,9)

                              \pspolygon[fillcolor=gray, linecolor=gray,
                              fillstyle=solid](3,9)(3,5)(7,5)(7,2)(9,2)(9,9)

                              \psline[linecolor=darkgray,linewidth=2pt](3,9)(3,5)(7,5)(7,2)(9,2)

                              \psline[linecolor=black,linewidth=1pt](1,9)(1,6)(2,6)(2,4)(4,4)(4,1)(9,1)

                               \rput(1,6){{\Large $\bullet$}}%
                               \rput(2,4){{\Large $\bullet$}}%
                               \rput(3,5){{\Large $\bullet$}}%
                               \rput(4,1){{\Large $\bullet$}}%
                               \rput(5,8){{\Large $\bullet$}}%
                               \rput(6,7){{\Large $\bullet$}}%
                               \rput(7,2){{\Large $\bullet$}}%
                               \rput(8,3){{\Large $\bullet$}}%

                          \end{pspicture}\\[0.25in]

                      (c) Shadowline $L_{2}(64518723)$.
                  \end{minipage}
                  &
                  \begin{minipage}[c]{1.325in}
                      \centering

                          \psset{xunit=0.125in,yunit=0.125in}

                          \begin{pspicture}(0,0)(9,9)

                              \psaxes[Dy=2]{->}(9,9)

                              \pspolygon[fillcolor=gray, linecolor=gray,
                              fillstyle=solid](5,9)(5,8)(6,8)(6,7)(8,7)(8,3)(9,3)(9,9)

                              \psline[linecolor=darkgray,linewidth=2pt](5,9)(5,8)(6,8)(6,7)(8,7)(8,3)(9,3)

                              \psline[linecolor=black,linewidth=1pt](3,9)(3,5)(7,5)(7,2)(9,2)

                              \psline[linecolor=black,linewidth=1pt](1,9)(1,6)(2,6)(2,4)(4,4)(4,1)(9,1)

                               \rput(1,6){{\Large $\bullet$}}%
                               \rput(2,4){{\Large $\bullet$}}%
                               \rput(3,5){{\Large $\bullet$}}%
                               \rput(4,1){{\Large $\bullet$}}%
                               \rput(5,8){{\Large $\bullet$}}%
                               \rput(6,7){{\Large $\bullet$}}%
                               \rput(7,2){{\Large $\bullet$}}%
                               \rput(8,3){{\Large $\bullet$}}%

                          \end{pspicture}\\[0.25in]

                      (d) Shadowline $L_{3}(64518723)$.
                  \end{minipage}

              \end{tabular}\\

              \caption{Examples of Northeast Shadow and Shadowline
              Constructions}

          \end{figure}
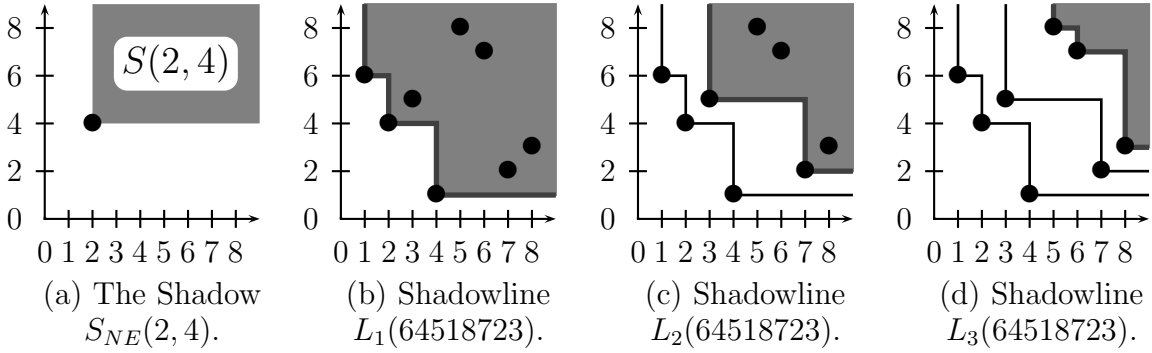

        At the same time, it is also clear that in general there will
        be many permutations $\sigma, \tau \in \mathfrak{S}_{n}$ for
        which $R(\sigma) = R(\tau)$.  In Section
        \ref{sec:PileConfigurations:PSEquivalence} below we
        characterize when two permutations have the same pile
        configuration, and we will denote this equivalence relation by
        $\sigma \stackrel{PS}{\sim} \tau$.  Moreover, we will also see
        that the reverse patience word $RPW(R(\sigma))$ is the most
        natural representative for the equivalence class generated by
        a given permutation $\sigma$.\bigskip

        We close this section by giving an alternate characterization
        for pile configurations in terms of the \emph{shadow diagram}
        construction that Viennot \cite{refViennot1977} introduced in
        the context of studying the RSK Algorithm for permutations.

        \begin{defn}\label{defn:Shadow}
            Given a lattice point $(m, n) \in \mathbb{Z}^{2}$, we
            define the \emph{(northeast) shadow} of $(m, n)$ to be the
            quarter space $S(m, n) = \{ (x, y) \in \mathbb{R}^{2} \ |
            \ x \geq m, \ y \geq n\}$.
        \end{defn}

        \noindent See Figure \ref{fig:ShadowExample}(a) for an example
        of a point's shadow.\medskip

        The most important use of shadows is in building shadowlines:

        \begin{defn}\label{defn:Shadowline}
            Given lattice points $(m_{1}, n_{1}), (m_{2}, n_{2}),
            \ldots, (m_{k}, n_{k}) \in \mathbb{Z}^{2}$, we define
            their \emph{(northeast) shadowline} to be the boundary of
            the region formed from the union of the shadows $S(m_{1},
            n_{1}), S(m_{2}, n_{2}), \ldots, S(m_{k}, n_{k})$.
        \end{defn}

        In particular, we wish to associate to each permutation a
        certain collection of shadowlines (as illustrated in Figure
        \ref{fig:ShadowExample}(b)--(d)):

        \begin{defn}\label{defn:ShadowDiagram}
            Given a permutation $\sigma =
            \sigma_{1}\sigma_{2}\cdots\sigma_{n} \in
            \mathfrak{S}_{n}$, the \emph{(northeast) shadow diagram}
            $D(\sigma)$ of $\sigma$ consists of the shadowlines
            $L_{1}(\sigma), L_{2}(\sigma), \ldots, L_{k}(\sigma)$
            formed as follows:\medskip

            \begin{itemize}
                \item $L_{1}(\sigma)$ is the shadowline for the
                lattice points $\{(1, \sigma_{1}), (2, \sigma_{2}),
                \ldots, (n, \sigma_{n})\}$.\medskip

                \item While at least one of the points $(1,
                \sigma_{1}), (2, \sigma_{2}), \ldots, (n, \sigma_{n})$
                is not contained in the shadowlines $L_{1}(\sigma),
                L_{2}(\sigma), \ldots, L_{j}(\sigma)$, define
                $L_{j+1}(\sigma)$ to be the shadowline for the points
                $$\{(i, \sigma_{i}) \ | \ (i, \sigma_{i})
                \notin \bigcup^{j}_{k=1} L_{k}(\sigma)\}.$$
            \end{itemize}

        \end{defn}

        In other words, we define the shadow diagram $D(\sigma) = \{
        L_{1}(\sigma), L_{2}(\sigma), \ldots, L_{k}(\sigma) \}$
        inductively with $L_{1}(\sigma)$ the shadowline for the
        diagram $\{(1, \sigma_{1}), (2, \sigma_{2}), \ldots, (n,
        \sigma_{n})\}$ of the permutation $\sigma \in
        \mathfrak{S}_{n}$.  Then we ignore the points whose shadows
        were actually used in building $L_{1}(\sigma)$ and define
        $L_{2}(\sigma)$ to be the shadowline of the resulting subset
        of the permutation diagram.  We then build $L_{3}(\sigma)$ as
        the shadowline for the points not yet used in constructing
        both $L_{1}(\sigma)$ and $L_{2}(\sigma)$, and this process
        continues until all points in the permutation diagram are
        exhausted.\medskip

        One of the most basic properties of the shadow diagram for a
        permutation $\sigma$ is that it encodes the top row of the
        insertion tableau $P(\sigma)$ (resp.  recording tableau
        $Q(\sigma)$) as the smallest ordinates (resp.  smallest
        abscissae) of all points belonging to the constituent
        shadowlines $L_{1}(\sigma), L_{2}(\sigma), \ldots,
        L_{k}(\sigma)$.  (A proof of this can be found in Sagan
        \cite{refSagan2000}.)  In particular, this means that if
        $\sigma$ has pile configuration $R(\sigma) = \{r_{1}, r_{2},
        \ldots, r_{m}\}$, then $m = k$ since the number of piles is
        equal to the length of the top row of $P(\sigma)$ (as both are
        the length of the longest increasing subsequence of $\sigma$).
        We can say even more about the relationship between
        $D(\sigma)$ and $R(\sigma)$ when both are viewed in terms of
        left-to-right minima subsequences (a.k.a. basic subsequences
        or records):

        \begin{defn}\label{defn:LtoRminimaSubsequence}
            Let $\pi = \pi_{1} \pi_{2} \cdots \pi_{l}$ be a partial
            permutation on $[n] = \{1, 2, \ldots, n\}$.  Then the
            \emph{left-to-right minima subsequence} of $\pi$ consists
            of those components $\pi_{j}$ of $\pi$ such that $\pi_{j}
            = \min\{ \pi_{i} \ | \ 1 \leq i \leq j\}$.
        \end{defn}

        \noindent We then inductively define the left-to-right minima
        subsequences $s_{1}, s_{2}, \ldots, s_{k}$ of a permutation
        $\sigma$ by taking $s_{1}$ to be the left-to-right minima
        subsequence for $\sigma$ itself and then each subsequent
        subsequence $s_{i}$ to be the left-to-right minima subsequence
        for the partial permutation obtained by removing the elements
        of $s_{1}, s_{2}, \ldots, s_{i-1}$ from $\sigma$.  The
        subsequence $s_j$ is called the $j^{\mathrm{th}}$
        \emph{left-to-right minima subsequence} of $\sigma$.

        \begin{lem}\label{lem:ShadowDiagramPileCorrespondence}
            Suppose $\sigma \in \mathfrak{S}_{n}$ has shadow diagram
            $D(\sigma) = \{ L_{1}(\sigma), L_{2}(\sigma), \ldots,
            L_{k}(\sigma) \}$.  Then the ordinates of the southwest
            corners of $L_{j}$ are exactly the cards in the
            $j^{\mathrm{th}}$ pile $r_{j} \in R(\sigma)$ formed by
            applying Patience Sorting (Algorithm
            \ref{alg:MallowsPSprocedure}) to $\sigma$.  In other
            words, the $j^{\mathrm{th}}$ pile $r_j$ contains exactly
            the elements of the $j^{\mathrm{th}}$ left-to-right minima
            subsequence of $\sigma$.
        \end{lem}

        \begin{proof}
            The $i^{\mathrm{th}}$ basic subsequence $s_{i}$ of $\sigma$
            consists of those elements $\sigma_{t}$ that appear at the
            end of an increasing subsequence of length $i$ but not at
            the end of an increasing subsequence of length $i+1$.
            Thus, since each element added to a pile must be smaller
            than all other elements already in the pile, $s_{1} =
            r_{1}$.  It then follows similarly by induction that
            $s_{i} = r_{i}$ for $i = 2, \ldots, k$.

            The proof that the ordinates of the southwest corners of
            the shadowlines $L_{i}$ are also exactly the elements of
            the left-to-right minima subsequences $s_{i}$ is similar.
        \end{proof}

        Lemma \ref{lem:ShadowDiagramPileCorrespondence} gives a
        particularly nice correspondence between the piles formed
        under Patience Sorting and the shadowlines that constitute the
        shadow diagram of a permutation.  In particular, we have that
        forming $RPW(R(\sigma))$ essentially amounts to sorting
        $\sigma \in \mathfrak{S}_{n}$ into its left-to-right minima
        subsequences.

        \medskip

        We will rely heavily upon this correspondence in the sections
        below.

        \subsection{Permutations Having Equivalent Pile Configurations}
        \label{sec:PileConfigurations:PSEquivalence}

        In this section we characterize the following equivalence
        relation:

        \begin{defn}\label{defn:PSequivalence}
            Two permutations $\sigma, \tau \in \mathfrak{S}_{n}$ are
            said to be \emph{patience sorting equivalent}, written
            $\sigma \stackrel{PS}{\sim} \tau$, if they have the same
            pile configuration $R(\sigma) = R(\tau)$ under Algorithm
            \ref{alg:MallowsPSprocedure}.  We denote the equivalence
            class generated by $\sigma$ as $\widetilde{\sigma}$.
        \end{defn}

        \begin{figure}[t]\label{fig:PSequivalenceExamples}
              \centering

              \begin{tabular}{cc}

                  \begin{minipage}[c]{3in}
                      \centering
                      \begin{tabular}{c c c}

                          \psset{xunit=0.25in,yunit=0.25in}

                          \begin{pspicture}(0,0)(4,4)

                              \psaxes{->}(4,4)

                              \psline[linecolor=black,linewidth=1pt](1,4)(1,2)(3,2)(3,1)(4,1)

                              \psline[linecolor=black,linewidth=1pt](2,4)(2,3)(4,3)

                               \rput(1,2){{\Large $\bullet$}}%
                               \rput(2,3){{\Large $\bullet$}}%
                               \rput(3,1){{\Large $\bullet$}}%

                          \end{pspicture}

                           & \raisebox{0.35in}{$\stackrel{PS}{\sim}$ \
                           } &

                          \psset{xunit=0.25in,yunit=0.25in}

                          \begin{pspicture}(0,0)(4,4)

                              \psaxes{->}(4,4)

                              \psline[linecolor=black,linewidth=1pt](1,4)(1,2)(2,2)(2,1)(4,1)

                              \psline[linecolor=black,linewidth=1pt](3,4)(3,3)(4,3)

                               \rput(1,2){{\Large $\bullet$}}%
                               \rput(2,1){{\Large $\bullet$}}%
                               \rput(3,3){{\Large $\bullet$}}%

                          \end{pspicture}\\[0.35in]

                      \end{tabular}

                      (a) ``Stretching'' shadowlines effects \\ $231
                      \stackrel{PS}{\sim} 213$.  Thus, $\widetilde{231} =
                      \{231, 213\}$.
                  \end{minipage}
                  &
                  \begin{minipage}[c]{2.125in}
                      \centering

                          \psset{xunit=0.2in,yunit=0.2in}

                          \begin{pspicture}(0,0)(5,5)

                              \psaxes{->}(5,5)

                              \psline[linecolor=black,linewidth=1pt](1,5)(1,3)(2,3)(2,1)(5,1)

                              \psline[linecolor=black,linewidth=1pt](3,5)(3,4)(4,4)(4,2)(5,2)

                               \rput(1,3){{\Large $\bullet$}}%
                               \rput(2,1){{\Large $\bullet$}}%
                               \rput(3,4){{\Large $\bullet$}}%
                               \rput(4,2){{\Large $\bullet$}}%

                          \end{pspicture}\\[0.35in]

                          (b) No ``stretching'' can interchange ``4''
                          and ``2''.
                  \end{minipage}

                  \end{tabular}

              \caption{Examples of patience sorting equivalence and
              non-equivalence}

          \end{figure}
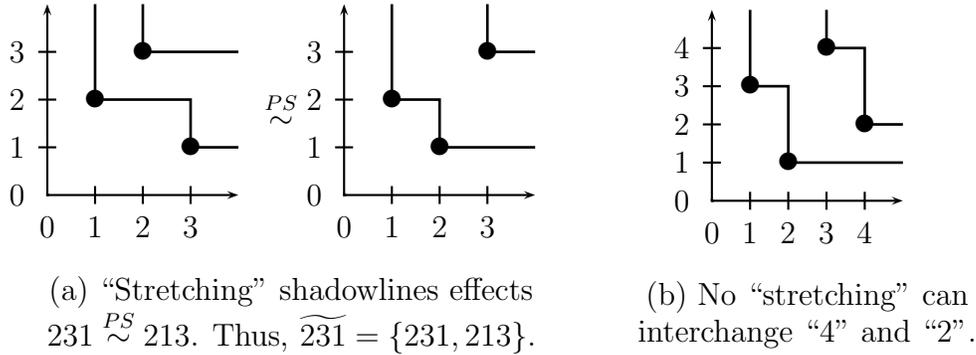

        By Lemma \ref{lem:ShadowDiagramPileCorrespondence} in Section
        \ref{sec:PileConfigurations:Characterization} above, the pile
        configurations $R(\sigma)$ and $R(\tau)$ correspond to certain
        shadow diagrams.  Thus, it should be intuitively clear that
        preserving a given pile configuration is equivalent to
        preserving the ordinates for the southwest corners of the
        shadowlines.  In particular, this means that we are limited to
        horizontally ``stretching'' shadowlines up to the point of not
        allowing them to cross as is illustrated in Figure
        \ref{fig:PSequivalenceExamples} and the following examples.

        \begin{exgr}\label{eg:PSequivalenceExample}
            The only non-singleton patience sorting equivalence class
            for $\mathfrak{S}_{3}$ consists of $\widetilde{231} = \{231,
            213\}$.  We illustrate $231 \stackrel{PS}{\sim} 213$ in
            Figure \ref{fig:PSequivalenceExamples}(a).
        \end{exgr}

        \noindent Notice that the actual values of the elements
        interchanged in Example \ref{eg:PSequivalenceExample} are
        immaterial so long as they have the same relative magnitudes
        as the literal values in the word $231$.  (I.e., they have to
        be order-isomorphic.)  Moreover, it should also be clear that
        any value greater than the element playing the role of ``1''
        can be inserted between the elements playing the roles of
        ``2'' and ``3'' without affecting the ability to interchange
        the ``1'' and ``3'' elements.  Problems with this interchange
        only start to arise when a value smaller than the element
        playing the role of ``1'' is inserted between the elements
        playing the roles of ``2'' and ``3''.  We can formally
        describe this idea using the language of generalized
        permutation patterns (as was recently defined in
        \cite{refBS2000}; cf.  \cite{refBona2004}).

        \begin{defn}\label{defn:PermutationPattern}
            Let $\sigma = \sigma_{1}\sigma_{2}\cdots\sigma_{n} \in
            \mathfrak{S}_{n}$ and $\tau \in \mathfrak{S}_{m}$ for $m
            \leq n$.  Then we say that $\sigma$ \emph{contains the
            (classical) pattern} $\tau$ if there exists a subsequence
            $\sigma_{i_{1}}, \sigma_{i_{2}}, \ldots, \sigma_{i_{m}}$
            of $\sigma$ (meaning $i_{1} < i_{2} < \cdots < i_{m}$)
            such that the word
            $\sigma_{i_{1}}\sigma_{i_{2}}\ldots\sigma_{i_{m}}$ is
            order-isomorphic to $\tau$.

            If $\sigma$ does not contain $\tau$, then we say that
            $\sigma$ \emph{avoids} the pattern $\tau$, and we denote
            by $S_{n}(\tau)$ the subset of the symmetric group
            $\mathfrak{S}_{n}$ that avoids $\tau$.
        \end{defn}

        Note that the elements in the subsequence $\sigma_{i_{1}},
        \sigma_{i_{2}}, \ldots, \sigma_{i_{m}}$ are not required to be
        contiguous in $\sigma$.  In a \emph{generalized pattern} one
        assumes that every element in the subsequence must be taken
        contiguously unless a dash is inserted in the pattern $\tau$
        between elements that are not required to be contiguous in
        $\sigma$.  (A generalized patterns with no dashes is sometimes
        called a \emph{segment} or a \emph{consecutive pattern}.)

         \begin{exgr}\label{eg:PSequivalenceNonExample}
         \begin{tabular}{c} $\phantom{foo}$ \end{tabular}

         \begin{enumerate}
            \item Notice that $2431$ contains exactly one instance of
            a $2\textrm{-}31$ pattern as the bold underlined
            subsequence
            $\mathbf{\underline{2}}4\mathbf{\underline{3}}\mathbf{\underline{1}}$.
            (Conversely,
            $\mathbf{\underline{2}}\mathbf{\underline{4}}3\mathbf{\underline{1}}$
            is an instance of $23\textrm{-}1$ but not of
            $2\textrm{-}31$.)  Moreover, it is clear that $2431
            \stackrel{PS}{\sim} 2413$.

            \item Even though $3142$ contains a $2\textrm{-}31$
            pattern (as the subsequence
            $\mathbf{\underline{3}}1\mathbf{\underline{4}}\mathbf{\underline{2}}$),
            we cannot interchange ``4'' and ``2'', and so $R(3142)
            \neq R(3124)$.  As illustrated in Figure
            \ref{fig:PSequivalenceExamples}(b), this is because ``4''
            and ``2'' are on the same shadowline.
         \end{enumerate}

        \end{exgr}

          \noindent We can now state our main result on patience
          sorting equivalence:

          \begin{thm}\label{thm:PSequivalence}
            Let $\sigma, \tau \in \mathfrak{S}_{n}$.  Then $\sigma$
            and $\tau$ have the same pile configuration $R(\sigma) =
            R(\tau)$ under Algorithm \ref{alg:MallowsPSprocedure} (so
            that $\sigma \stackrel{PS}{\sim} \tau$) if and only if
            there exists a sequence of $2\textrm{-}31$ to
            $2\textrm{-}13$ interchanges (with no $2\textrm{-}31$
            pattern contained in a $3\textrm{-}1\textrm{-}42$ pattern)
            that transform $\sigma$ into $\tau$.

            In other words, $\stackrel{PS}{\sim}$ is the transitive
            closure of such interchanges.
        \end{thm}

        \begin{proof}
            By Lemma \ref{lem:ShadowDiagramPileCorrespondence} it
            suffices to show that $2\textrm{-}31$ to $2\textrm{-}13$
            interchanges (with no $2\textrm{-}31$ pattern contained in
            a $3\textrm{-}1\textrm{-}42$ pattern), preserve the
            left-to-right minima subsequences $s_{1}, s_{2}, \ldots,
            s_{k}$ of $\sigma$.  This amounts to showing by induction
            that such interchanges suffice to transform $\sigma$ into
            $RPW(R(\sigma))$ via the sequence of pattern interchanges
            \[
            \sigma = \sigma^{(0)} \leadsto \sigma^{(1)} \leadsto
            \sigma^{(2)} \leadsto \cdots \leadsto \sigma^{(\ell)} =
            RPW(R(\sigma))
            \]
            where each $\sigma^{(i)} \stackrel{PS}{\sim} \sigma^{(i+1)}$.

            Let $(a,b,c)$ be a subsequence of $\sigma_i$ that is an
            instance of $2\text{-}31$ not contained in an instance of
            $3\text{-}1\text{-}42$.  Then $c<a<b$, and there is no $d$
            between $a$ and $b$ in $\sigma$ such that $d<c$.  Clearly,
            we do not lose any increasing subsequences of
            $\sigma^{(i)}$ by interchanging $b$ and $c$.  Moreover,
            the only new increasing subsequences $\sigma^{(i+1)}$
            created by this interchange are those that end with the
            subsequence $(c,b)$.  If such an increasing subsequence
            contains other terms, say $d_1<d_2<\dots<d_m<c<b$, then
            $d_m$ must be to the left of $a$ in $\sigma$.  But then
            $\sigma^{(i)}$ contains an increasing sequence
            $d_1<d_2<\dots<d_m<a<b$ with the same number of terms.
            Hence, interchanging $b$ and $c$ does not create any
            longer increasing subsequences with the same final term,
            so $\sigma^{(i)}\stackrel{PS}{\sim}\sigma^{(i+1)}$.

            Let $r_{11}>r_{12}>\dots>r_{1s_1}$ be the subsequence of
            left-to-right minima of $\sigma$; from the proof of Lemma
            \ref{lem:ShadowDiagramPileCorrespondence}, these are the
            entries that form the leftmost pile $r_1$ of $R(\sigma)$.
            Now suppose that, for some $j<s_1$, there is an entry
            between $r_{1 j}$ and $r_{1, j+1}$, and let $b$ be the
            entry immediately preceding $r_{1, j+1}$.  Then $b > r_{1
            j}$ so that $(r_{1 j},b,r_{1, j+1})$ is an instance of
            $2\text{-}31$.  On the other hand, $r_{1, j+1}$ is the
            leftmost entry of $\sigma$ that is less than $r_{1 j}$, so
            no entry $d<r_{1, j+1}$ may occur between $r_{1 j}$ and
            $b$.  Hence, $(r_{1 j}, b, r_{1, j+1})$ is not an instance
            of $3\text{-}1\text{-}42$, and so interchanging $r_{1,
            j+1}$ and $b$ will not change the pile configuration.  We
            may repeat this until there are no entries between
            consecutive left-to-right minima of $\sigma$.  We thus
            obtain $\sigma'=r_{11}r_{12}\dots r_{1s_1}\sigma''$, where
            $\sigma''$ is obtained by deleting
            $r_{11},r_{12},\dots,r_{1s_1}$ from $\sigma$.  Since no
            instance of $2\text{-}31$ may start with $r_{1 s_1}$ and
            since any instance $(r_{1i}, a, b)$ of $2\text{-}31$ is
            part of an instance $(r_{1i}, r_{1 s_1}, a, b)$ of
            $3\text{-}1\text{-}42$, no further interchanges will
            involve any $r_{1i}$.  Thus, by induction, we can now
            apply the same procedure to $\sigma''$, etc., to
            ultimately obtain $r_{1 1} r_{1 2} \dots r_{1 s_{1}} \,
            r_{2 1} r_{2 2} \dots r_{2 s_{2}} \; \dots \; r_{k 1} r_{k
            2} \dots r_{k s_{k}}=RPW(R(\sigma))$.
        \end{proof}

        \begin{rem}
            It follows from Theorem \ref{thm:PSequivalence} that
            Examples \ref{eg:PSequivalenceExample} and
            \ref{eg:PSequivalenceNonExample}(2) sufficiently
            characterize when two permutations yield the same pile
            configurations under Patience Sorting.  However, it is
            worth pointing out that these examples also begin to
            illustrate how one can build an infinite sequence of
            generalized permutation patterns (all of them containing
            either $2\textrm{-}13$ or $2\textrm{-}31$) with the
            following property: an interchange of the pattern
            $2\textrm{-}13$ with the pattern $2\textrm{-}31$ is
            allowed within an odd-length pattern in this sequence
            unless the elements used to form the odd-length pattern
            can also be used as part of a longer even-length pattern
            in this sequence.
        \end{rem}

        \begin{exgr}
            Even though the permutation $34152$ contains a
            $3\textrm{-}1\textrm{-}42$ pattern in the suffix ``4152",
            one can still directly interchange the ``5" and the ``2"
            because of the ``3" prefix (or via the following
            sequence of interchanges: $34152 \leadsto 31452 \leadsto
            31425 \leadsto 34125$).
        \end{exgr}

    \section{Bijectively Extending Patience Sorting to ``Stable Pairs'' of Pile Configurations}
    \label{sec:ExtendingPS}

        \subsection{The Extended Patience Sorting Algorithm}
        \label{sec:ExtendingPS:Algorithm}

        Recall from Section \ref{sec:Introduction} that Patience
        Sorting (Algorithm \ref{alg:MallowsPSprocedure}) can be viewed
        as an iterated, non-recursive form of the Schensted Insertion
        Algorithm for inserting a value into the top row of a Young
        tableau.  In this section we extend the Patience Sorting
        construction so that it becomes a full non-recursive analog of
        the RSK Algorithm for permutations.  In particular, we mimic
        the RSK recording tableau construction so that ``recording
        piles'' are formed while assembling the usual pile
        configuration under Patience Sorting (which by analogy to RSK
        we will similarly now call ``insertion piles'').

        \begin{alg}[Extended Patience Sorting
        Algorithm]\label{alg:ExtendedPSalgorithm} Given a shuffled
        deck of cards $\sigma = c_{1} c_{2} \cdots c_{n}$, inductively
        build \emph{insertion piles} $R = R(\sigma) = \{r_{1}, r_{2},
        \ldots, r_{m}\}$ and \emph{recording piles} $S = S(\sigma) =
        \{s_{1}, s_{2}, \ldots, s_{m}\}$ as follows:\smallskip

            \begin{itemize}
                \item Place the first card $c_{1}$ from the deck into
                a pile $r_{1}$ by itself and set $s_{1} =
                \{1\}$.\smallskip

                \item For each remaining card $c_{i}$ ($i = 2, \ldots,
                n$), consider the cards $d_{1}, d_{2}, \ldots, d_{k}$
                atop the piles $r_{1}, r_{2}, \ldots, r_{k}$ that have
                already been formed.\smallskip

                    \begin{itemize}
                        \item If $c_{i} > \max\{d_{1}, d_{2}, \ldots,
                        d_{k}\}$, then put $c_{i}$ into a new pile
                        $r_{k+1}$ by itself and set $s_{k+1} =
                        \{i\}$.\smallskip

                        \item Otherwise, find the left-most card
                        $d_{j}$ that is larger than $c_{i}$ and put
                        the card $c_{i}$ atop pile $r_{j}$ while
                        simultaneously putting $i$ at the bottom of
                        pile $s_{j}$.
                    \end{itemize}
            \end{itemize}
        \end{alg}

        We call the pile configuration pairs that result from
        Algorithm \ref{alg:ExtendedPSalgorithm} \emph{stable pairs}
        and give a characterization for them in Section
        \ref{sec:ExtendingPS:CharacterizingStablePairs} below.  Note
        that the pile configurations that comprise a resulting stable
        pair must have the same ``shape'', which we define as follows:

        \begin{defn}
            Given a pile configuration $R=\{r_{1}, r_{2}, \ldots,
            r_{m}\}$ on $n$ cards, we call the composition $\gamma =
            (|r_{1}|, |r_{2}|, \ldots, |r_{m}|)$ of $n$ the
            \emph{shape} of $R$ and denote this by $\mathrm{sh}(R) =
            \gamma \models n$.
        \end{defn}

        \begin{exgr}\label{eg:ExtendedPSexample}
            Let $\sigma = 6 4 5 1 8 7 2 3 \in \mathfrak{S}_{8}$.  Then
            according to Algorithm \ref{alg:ExtendedPSalgorithm} we
            simultaneously form the following pile configurations with
            shape $\mathrm{sh}(R(\sigma)) = \mathrm{sh}(S(\sigma)) =
            (3, 2, 3) \models 8$.\\

            \begin{tabular}{l p{56pt} p{72pt} l p{56pt} l}
                \begin{minipage}[c]{48pt}$\phantom{foo}$\end{minipage}
                &
                \begin{minipage}[c]{48pt}insertion piles\end{minipage}
                &
                \begin{minipage}[c]{48pt}recording piles\end{minipage}
                &
                \begin{minipage}[c]{48pt}$\phantom{foo}$\end{minipage}
                &
                \begin{minipage}[c]{48pt}insertion piles\end{minipage}
                &
                \begin{minipage}[c]{48pt}recording piles\end{minipage}
             \end{tabular}\\

            \begin{tabular}{l p{56pt} p{72pt} l p{56pt} l}
                \begin{minipage}[c]{48pt}
                    Form a new pile with \textbf{6}:
                \end{minipage}
                &
                \begin{tabular}{l l l}
                                      & & \\
                                      & & \\
                       \textbf{6} & &
                \end{tabular}
                &
                \begin{tabular}{l l l}
                                     & & \\
                                     & & \\
                      \textbf{1} & &
                \end{tabular}
                &
                \begin{minipage}[c]{48pt}
                    Then play the \textbf{4} on it:
                \end{minipage}
                &
                \begin{tabular}{l l l}
                                      & & \\
                       \textbf{4} & & \\
                                  6 & &
                \end{tabular}
                &
                \begin{tabular}{l l l}
                                   & & \\
                               1 & & \\
                    \textbf{2} & &
                \end{tabular}
            \end{tabular}\\ \\

            \begin{tabular}{l p{56pt} p{72pt} l p{56pt} l}
                \begin{minipage}[c]{48pt}
                    Form a new pile with \textbf{5}:
                \end{minipage}
                &
                \begin{tabular}{l l l}
                       & & \\
                    4 & & \\
                    6 & \textbf{5} &
                \end{tabular}
                &
                \begin{tabular}{l l l}
                      & & \\
                   1 & & \\
                   2 & \textbf{3} &
                \end{tabular}
                &
                \begin{minipage}[c]{48pt}
                    Add the \textbf{1} to left pile:
                \end{minipage}
                &
                \begin{tabular}{l l l}
                    \textbf{1} & & \\
                               4 & & \\
                               6 & 5 &
                \end{tabular}
                &
                \begin{tabular}{l l l}
                               1 & & \\
                               2 & & \\
                    \textbf{4} & 3 &
                \end{tabular}
            \end{tabular}\\ \\

            \begin{tabular}{l p{56pt} p{72pt} l p{56pt} l}
                \begin{minipage}[c]{48pt}
                    Form a new pile with \textbf{8}:
                \end{minipage}
                &
                \begin{tabular}{l l l}
                    1 & & \\
                    4 & & \\
                    6 & 5 & \textbf{8}
                \end{tabular}
                &
                \begin{tabular}{l l l}
                   1 & & \\
                   2 & & \\
                   4 & 3 & \textbf{5}
                \end{tabular}
                &
                \begin{minipage}[c]{48pt}
                    Then play the \textbf{7} on it:
                \end{minipage}
                &
                \begin{tabular}{l l l}
                    1 & & \\
                    4 & & \textbf{7} \\
                    6 & 5 & 8
                \end{tabular}
                &
                \begin{tabular}{l l l}
                    1 & & \\
                    2 & & 5 \\
                    4 & 3 & \textbf{6}
                \end{tabular}
            \end{tabular}\\ \\

            \begin{tabular}{l p{56pt} p{72pt} l p{56pt} l}
                \begin{minipage}[c]{48pt}
                    Add the \textbf{2} to middle pile:
                \end{minipage}
                &
                \begin{tabular}{l l l}
                    1 & & \\
                    4 & \textbf{2} & 7 \\
                    6 & 5 & 8
                \end{tabular}
                &
                \begin{tabular}{l l l}
                   1 & & \\
                   2 & 3 & 5 \\
                   4 & \textbf{7} & 6
                \end{tabular}
                &
                \begin{minipage}[c]{48pt}
                    Add the \textbf{3} to right pile:
                \end{minipage}
                &
                \begin{tabular}{l l l}
                    1 & & \textbf{3} \\
                    4 & 2 & 7 \\
                    6 & 5 & 8
                \end{tabular}
                &
                \begin{tabular}{l l l}
                   1 & & 5 \\
                   2 & 3 & 6 \\
                   4 & 7 & \textbf{8}
                \end{tabular}
            \end{tabular}\\ \\

        \end{exgr}

        The idea behind Algorithm \ref{alg:ExtendedPSalgorithm} is
        that we are using the recording piles $S(\sigma)$ to
        implicitly label the order in which the elements of the
        permutation $\sigma$ are added to the insertion piles
        $R(\sigma)$.  It is clear that this information then allows us
        to uniquely reconstruct $\sigma$ by reversing the order in
        which the cards were played.  However, even though reversing
        the Extended Patience Sorting Algorithm is much easier than
        reversing the RSK Algorithm through recursive ``reverse row
        bumping'', the trade-off is that the stable pairs that result
        from the former are not independent whereas the tableau pairs
        generated by RSK are completely independent (up to shape).

        That $S(\sigma) = \{s_{1}, s_{2}, \ldots, s_{m}\}$ records the
        order of the cards being added to the insertion piles is made
        clear if we instead add cards to the tops of new piles
        $s'_{j}$ in Algorithm \ref{alg:ExtendedPSalgorithm} rather
        than to the bottoms of the piles $s_{j}$.  This yields
        modified recording piles $S'(\sigma)$ from which each original
        recording pile $s_{j} \in S(\sigma)$ can be recovered by
        simply reflecting the corresponding pile $s'_{j}$ vertically.

        \begin{exgr}
            As in Example \ref{eg:ExtendedPSexample} above, let
            $\sigma = 6 4 5 1 8 7 2 3 \in \mathfrak{S}_{8}$.  Then
            $R(\sigma)$ is formed as before and

            \begin{center}
                \begin{tabular}{r c c c l}
                    $S'(\sigma) \ \ =$ &
                    \begin{tabular}{l l l}
                        4 & & 8 \\
                        2 & 7 & 6 \\
                        1 & 3 & 5
                    \end{tabular}
                    & $\stackrel{reflect}{\dashrightarrow}$ &
                    \begin{tabular}{l l l}
                       1 & & 5 \\
                       2 & 3 & 6 \\
                       4 & 7 & 8
                    \end{tabular}
                    & $= \ \ S(\sigma)$
                \end{tabular}
            \end{center}

        \end{exgr}

        We are now in a position to prove that the Extended Patience
        Sorting Algorithm has the same form of symmetry as the RSK
        Algorithm has for permutations.

        \begin{prop}\label{prop:TwiddlePilesInverseResult}
            Let $(R(\sigma), S(\sigma))$ be the insertion and
            recording piles, respectively, formed by applying
            Algorithm \ref{alg:ExtendedPSalgorithm} to $\sigma \in
            \mathfrak{S}_{n}$.  Then reversing Algorithm
            \ref{alg:ExtendedPSalgorithm} for the pair
            $(S(\sigma),R(\sigma))$ yields the inverse permutation
            $\sigma^{-1}$.
        \end{prop}

        \begin{proof}
            Construct $S'(\sigma)$ from $S(\sigma)$ as discussed
            above, and form the $n$ ordered pairs $(r_{i j}, s'_{i
            j})$ where $i$ indexes the individual piles and $j$ the
            cards in the $i^{\mathrm{th}}$ piles.  Then these $n$
            points correspond to the diagram of a permutation $\tau
            \in \mathfrak{S}_{n}$.  However, since reflecting these
            points through the line $y = x$ yields the diagram for
            $\sigma$, it follows that $\tau = \sigma^{-1}$.
        \end{proof}

        Proposition \ref{prop:TwiddlePilesInverseResult} suggests that
        Algorithm~\ref{alg:ExtendedPSalgorithm} is the right
        generalization of Algorithm~\ref{alg:MallowsPSprocedure} since
        we obtain the same symmetry property as for RSK. At the same
        time, though, since there are many more possible pile
        configurations than standard Young tableau (as we'll show in
        Section \ref{sec:EnumeratingSn3142} below), not every ordered
        pair of pile configurations with the same shape will result
        from Algorithm \ref{alg:ExtendedPSalgorithm}.  Thus, it is
        necessary to first characterize the ``stable pairs'' that
        result from applying Extended Patience Sorting to a
        permutation.  We do this in Section
        \ref{sec:ExtendingPS:CharacterizingStablePairs}.

        \subsection{Characterizing ``Stable Pairs'' of Pile Configurations and Pile Configurations for Involutions}
        \label{sec:ExtendingPS:CharacterizingStablePairs}

        Based upon Proposition \ref{prop:TwiddlePilesInverseResult}
        above, there is a bijection between involutions and certain
        pile configurations.  We will describe this bijection as a
        corollary to a more general characterization for the ``stable
        pairs'' of pile configurations that can result from applying
        the Extended Patience Sorting Algorithm to a permutation.

        \medskip

        \noindent The following example, though very small,
        illustrates the most generic behavior that must be avoided in
        constructing stable pairs.  As in Section
        \ref{sec:ExtendingPS:Algorithm} above, we denote by $S'$ the
        ``reverse pile configuration'' of $S$ (which has all piles
        listed in reverse order).

        \begin{exgr}\label{eg:SmallBadPilesExample}

            Even though the pile configuration $R = \{\{3 > 1\},
            \{2\}\}$ cannot result as the insertion piles given by an
            involution under the Extended Patience Sorting Algorithm,
            we can still try to look at the shadow diagram for the
            pre-image of the pair $(R, R)$ under Algorithm
            \ref{alg:ExtendedPSalgorithm}:

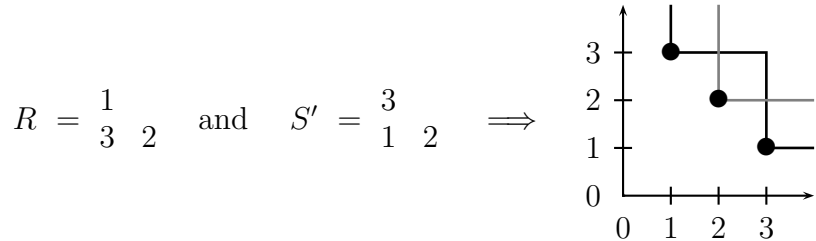
            \begin{center}

                \begin{tabular}{lr}

                    \raisebox{.35in}{
                            \begin{minipage}[c]{3in}

                                    $R \ = \
                                      \begin{matrix}
                                          1 & \\
                                          3 & 2
                                      \end{matrix}$
                                      \quad \text{and} \quad
                                      $S' \ = \
                                      \begin{matrix}
                                          3 & \\
                                          1 & 2
                                      \end{matrix}
                                      \quad \implies $

                            \end{minipage}
                    } &

                    \psset{xunit=0.25in,yunit=0.25in}

                    \begin{pspicture}(0,0)(4,4)

                        \psaxes{->}(4,4)

                        \psline[linecolor=black,linewidth=1pt](1,4)(1,3)(3,3)(3,1)(4,1)

                        \psline[linecolor=gray,linewidth=1pt](2,4)(2,2)(4,2)

                         \rput(1,3){{\Large $\bullet$}}%
                         \rput(2,2){{\Large $\bullet$}}%
                         \rput(3,1){{\Large $\bullet$}}%

                    \end{pspicture}\\[0.25in]

                \end{tabular}

            \end{center}

            \noindent Note that there are two competing constructions
            here.  On the one hand, we have the diagram $\{(1, 3), (2, 2),
            (3, 1)\}$ of a permutation given by the entries in the pile
            configurations.  (In particular, the values in $R$ specify the
            ordinates and the values in the corresponding boxes of $S'$
            the abscissae.)  On the other hand, the piles in $R$ also
            specify shadowlines with respect to this permutation diagram.
            Here the pair $(R, S)$ of pile configurations is ``unstable''
            because their combination yields crossing shadowlines---which
            is clearly not allowed.\\
            \indent Similar considerations lead to avoiding the crossings of
            the form\\

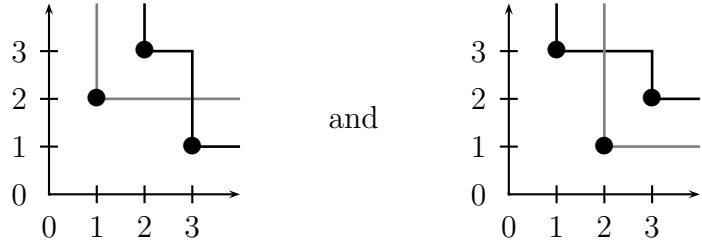
            \begin{center}

                \begin{tabular}{rlr}

                    \psset{xunit=0.25in,yunit=0.25in}

                    \begin{pspicture}(0,0)(4,4)

                        \psaxes{->}(4,4)

                        \psline[linecolor=gray,linewidth=1pt](1,4)(1,2)(4,2)

                        \psline[linecolor=black,linewidth=1pt](2,4)(2,3)(3,3)(3,1)(4,1)

                         \rput(1,2){{\Large $\bullet$}}%
                         \rput(2,3){{\Large $\bullet$}}%
                         \rput(3,1){{\Large $\bullet$}}%

                    \end{pspicture}

                    &

                    \raisebox{.35in}{\qquad and \qquad\quad}

                    &

                    \psset{xunit=0.25in,yunit=0.25in}

                    \begin{pspicture}(0,0)(4,4)

                        \psaxes{->}(4,4)

                        \psline[linecolor=black,linewidth=1pt](1,4)(1,3)(3,3)(3,2)(4,2)

                        \psline[linecolor=gray,linewidth=1pt](2,4)(2,1)(4,1)

                         \rput(1,3){{\Large $\bullet$}}%
                         \rput(2,1){{\Large $\bullet$}}%
                         \rput(3,2){{\Large $\bullet$}}%

                    \end{pspicture}\\[0.25in]

                \end{tabular}

            \end{center}

            \noindent Note also that these latter two crossings can
            also be used together to build something like the first
            crossing but with ``extra'' elements on the boundary of
            the polygon formed:

            \begin{center}

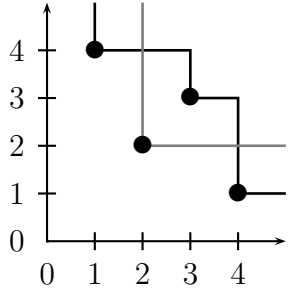
                \psset{xunit=0.25in,yunit=0.25in}

                \begin{pspicture}(0,0)(5,5)

                    \psaxes{->}(5,5)

                    \psline[linecolor=black,linewidth=1pt](1,5)(1,4)(3,4)(3,3)(4,3)(4,1)(5,1)

                    \psline[linecolor=gray,linewidth=1pt](2,5)(2,2)(5,2)

                     \rput(1,4){{\Large $\bullet$}}%
                     \rput(2,2){{\Large $\bullet$}}%
                     \rput(3,3){{\Large $\bullet$}}%
                     \rput(4,1){{\Large $\bullet$}}%

                \end{pspicture}\\[0.4in]

            \end{center}

        \end{exgr}

        We can now make the following important definitions:

        \begin{defn}
            Given a composition $\gamma$ of $n$ (denoted $\gamma
            \models n$), we define $\mathfrak{P}_{\gamma}(n)$ to be
            the set of all pile configurations $R$ having shape
            $\mathrm{sh}(R) = \gamma$ and put

            \begin{displaymath}
                \mathfrak{P}(n) = \bigcup_{\gamma
                \ \models \ n}\mathfrak{P}_{\gamma}(n).
            \end{displaymath}

        \end{defn}

        \begin{defn}\label{defn:StablePairs}
            Define the set $\Sigma(n) \subset \mathfrak{P}(n) \times
            \mathfrak{P}(n)$ to consist of all ordered pairs $(R, S)$
            with $\mathrm{sh}(R) = \mathrm{sh}(S)$ such that the
            pair $(RPW(R),RPW(S'))$ avoids simultaneous occurrences
            of pairs of patterns $(31\textrm{-}2,13\textrm{-}2)$,
            $(31\textrm{-}2,32\textrm{-}1)$ and $(32\textrm{-}1,13\textrm{-}2)$
            at the same positions in $RPW(R)$ and $RPW(S')$.
        \end{defn}

        In other words, if $RPW(R)$ contains an occurrence of the
        first pattern in any of the above pairs, then $RPW(S')$
        cannot contain an occurrence at the same positions of the
        second pattern in the same pair, and vice versa. Note that
        Definition \ref{defn:StablePairs}
        characterizes ``stable pairs'' of pile configurations $(R, S)$
        by forcing $R$ and $S$ to avoid certain sub-pile pattern
        pairs.  As in Example \ref{eg:SmallBadPilesExample}, we are
        characterizing when the induced shadowlines cross.

        \begin{thm}\label{thm:ExtendedPSbijection}
            Extended Patience Sorting (Algorithm
            \ref{alg:ExtendedPSalgorithm}) gives a bijection between
            the symmetric group $\mathfrak{S}_{n}$ and the ``stable
            pairs'' set $\Sigma(n)$ given in Definition
            \ref{defn:StablePairs} above.
        \end{thm}

        \begin{proof}
            We will prove that for any stable pair $(R,S)\in\Sigma(n)$
            and any permutation $\sigma\in \mathfrak{S}_{n}$,
            \[
            \sigma=\begin{pmatrix}RPW(S')\\RPW(R)\end{pmatrix}
            \; \text{(in the two-line notation)} \;
            \iff
            (R,S)=(R(\sigma),S(\sigma)).
            \]
            Clearly, if $(R,S)=(R(\sigma),S(\sigma))$ for some
            $\sigma\in \mathfrak{S}_{n}$, then
            $\sigma=\begin{pmatrix}RPW(S')\\RPW(R)\end{pmatrix}$, so
            we only need to prove that $(R,S)\in\Sigma(n)$.  Suppose
            that $(R,S)\notin\Sigma(n)$; then $RPW(R)$ and $RPW(S')$
            contain instances of one of the three forbidden pairs.

            Suppose $RPW(R)$ contains an occurrence $(r_3,r_1,r_2)$ of
            $31\text{-}2$, and $RPW(S')$ contains an occurrence
            $(s'_1,s'_3,s'_2)$ of $13\text{-}2$ at the
            same positions.  Since $r_3>r_1$ and $r_3,r_1$ are
            consecutive entries in $RPW(R)$, it follows that $r_3$ and
            $r_1$ must be in the same column $c_i(R)$ of $R$ (in fact,
            $r_1$ is immediately on top of $r_3$).  Since $r_1 < r_2$
            and $r_2$ is to the right of $r_1$ in $R$, it follows that
            the column $c_j(R)$ of $R$ containing $r_2$ must be to the
            right of the column containing $r_1$ atop $r_3$.
            Therefore, $s'_2$ must also be in a column $c_i(S')$ of
            $S'$ to the right of the column $c_j(S')$ containing
            $s'_3$ on top of $s'_1$.

            Consider the subpermutation $\tau$ of $\sigma$ formed by
            deleting entries of $RPW(R)$ and $RPW(S')$ that are not in
            these two columns.  Alternatively, let $R_\ast$ and
            $S'_\ast$ consist only of the columns $(c_i(R), c_j(R))$
            of $R$ and $(c_i(S'),c_j(S'))$ of $S'$, respectively.
            Then
            \[
            \tau=\begin{pmatrix}RPW(S'_\ast)\\RPW(R_\ast)\end{pmatrix}.
            \]
            Note that the values $r_3$ and $r_1$ in $c_i(S')$ are
            consecutive left-to-right minima of $\tau$, whereas $r_2$
            is not a left-to-right minimum of $\tau$.  Since
            $r_1<r_2<r_3$, it follows that $r_2$ cannot occur between
            $r_1$ and $r_3$ in $\tau$.  However, since
            $\begin{pmatrix}s'_1 & s'_3 & s'_2\\r_3 & r_1 &
            r_2\end{pmatrix}$ is a subpermutation of $\tau$ and
            $s'_1<s'_2<s'_3$, it follows that $r_2$ does occur between
            $r_1$ and $r_3$, a contradiction.

            A similar argument applies to the other two forbidden
            pairs, $(31\textrm{-}2,32\textrm{-}1)$ and
            $(32\textrm{-}1,13\textrm{-}2)$. The resulting contradiction
            in each case implies that we must have $(R,S)\in\Sigma(n)$.

            Conversely, if $(R,S)\in\Sigma(n)$, then set
            $\sigma=\begin{pmatrix}RPW(S')\\RPW(R)\end{pmatrix}$.  We
            must show that $(R,S)=(R(\sigma),S(\sigma))$. This can
            be done by showing that $c_1(R)$ is the sequence of
            left-to-right minima of $\sigma$ and $c_1(S')$ is the
            sequence of their positions, and then proceeding by
            induction on the length of $\sigma$.

            We know that $c_1(R)$ is decreasing and $c_1(S')$ is
            increasing, with both columns of the same length.
            Moreover, $c_1(R)$ is the sequence of left-to-right
            minima of $RPW(R)$, and $c_1(S')$ is increasing, so
            the values of $c_1(R)$ are in decreasing order in
            $\sigma$. Note that the first term of $c_1(R)$ is
            the leftmost value in $\sigma$, and the last term of
            $c_1(R)$ is the least value of $\sigma$.

            If every value in the first column of $R$ is less than
            every value not in the first column of $R$, and every
            value in the first column of $S'$ is less than
            every value not in the first column of $S'$,
            then we are done. Now suppose this is not so. Then there exists
            a point $(s',r)$ with $s'\in (RPW(S')\setminus c_1(S')$ and
            $r\in RPW(R)\setminus c_1(R))$
            such that there are points $(s'_1,r_1)$ and $(s'_2,r_2)$,
            where $r_1,r_2\in c_1(R)$ and $s'_1,s'_2\in c_1(S'))$ are
            consecutive terms occupying the same positions in the
            first columns of $R$ and $S'$, for which either $r_1>r>r_2$ or
            $s'_1<s'<s'_2$. None of the three forbidden pairs may
            occur at positions $(s'_1,r_1),(s'_2,r_2),(s',r)$, in
            other words, given that one of $r_1>r>r_2$ or
            $s'_1<s'<s'_2$ must hold, we must avoid the following
            cases:
            \[
            \begin{split}
            r_1>r>r_2 \ \text{and} \ s'<s'_1, &\qquad \text{i.e. an occurrence of
            $(31\textrm{-}2,32\textrm{-}1)$,}\\
            r_1>r>r_2 \ \text{and} \ s'_1<s'<s'_2, &\qquad \text{i.e. an occurrence of
            $(31\textrm{-}2,13\textrm{-}2)$,}\\
            r<r_1 \ \text{and} \ s'_1<s'<s'_2, &\qquad \text{i.e. an occurrence of
            $(32\textrm{-}1,13\textrm{-}2)$.}\\
            \end{split}
            \]
            Hence, $r_1>r>r_2$ implies $s'>s'_2$, while
            $s'_1<s'<s'_2$ implies $r>r_1$. Thus, $r_2$ is the first
            value of $\sigma$ to the right of $r_1$ that is less than
            $r_1$, and all values between $r_1$ and $r_2$ are to the
            right of $r_2$. Therefore, if $r_1$ is a left-to-right
            minimum of $\sigma$, then $r_2$ is the next left-to-right
            minimum of $\sigma$. Since the first term of $c_1(R)$ is
            a left-to-right minimum of $\sigma$, it follows by
            induction that $c_1(R)$ is the sequence of left-to-right
            minima of $\sigma$, which occur at positions in $c_1(S')$.
            This finishes the proof.
        \end{proof}

        \begin{exgr}\label{eg:GoodPilesExample}
        The pair of piles
        \[
        (R,S)=
        \left(
        \begin{tabular}{c c c}
            1 &   & 3 \\
            4 & 2 & 7 \\
            6 & 5 & 8
        \end{tabular}
        \; , \;
        \begin{tabular}{c c c}
            1 &   & 5 \\
            2 & 3 & 6 \\
            4 & 7 & 8
        \end{tabular}
        \right)
        \in\Sigma(8)
        \]
        corresponds to the permutation
        \[
        \sigma=\begin{pmatrix}RPW(S')\\RPW(R)\end{pmatrix}=
        \begin{pmatrix}
        1&2&4&3&7&5&6&8\\
        6&4&1&5&2&8&7&3
        \end{pmatrix}
        =
        64518723\in \mathfrak{S}_8.
        \]
        \end{exgr}

        Based upon the characterization of stable pairs given in
        Theorem \ref{thm:ExtendedPSbijection} and the
        Sch\"{u}t\-zen\-ber\-ger-type Symmetry Property proven in
        Proposition \ref{prop:TwiddlePilesInverseResult}, we can
        immediately describe a bijection between involutions and
        certain pile configurations.  In particular, these pile
        configurations must avoid simultaneously containing the
        symmetric sub-pile patterns corresponding to the patterns
        given in Definition \ref{defn:StablePairs}.

        This corresponds to the reverse patience word for a pile
        configuration simultaneously avoiding a symmetric pair of the
        generalized patterns $31\textrm{-}2$ and $32\textrm{-}1$.  As
        such it is interesting to compare this construction to two
        results recently obtained by Claesson and Mansour
        \cite{refCM}:

        \begin{enumerate}
            \item The size of $S_{n}(3\textrm{-}12, 3\textrm{-}21)$ is
            equal to the number of involutions $|\mathfrak{I}_{n}|$ in
            $\mathfrak{S}_{n}$.\smallskip

            \item The size of $S_{n}(31\textrm{-}2, 32\textrm{-}1)$ is
            $2^{n-1}$.
        \end{enumerate}

        \noindent The first result suggests that there should be a way
        to relate the result in Theorem~\ref{thm:ExtendedPSbijection}
        to simultaneous avoidance of the very similar patterns
        $3\textrm{-}12$ and $3\textrm{-}21$.  The second result
        suggests that restricting to complete avoidance of all
        simultaneous occurrences of $31\textrm{-}2$ and
        $32\textrm{-}1$ will yield a natural bijection between
        $S_{n}(31\textrm{-}2, 32\textrm{-}1)$ and a subset
        $\mathfrak{N} \subset \mathfrak{P}(n)$ such that $\mathfrak{N}
        \cap \mathfrak{P}_{\gamma}(n)$ contains exactly one pile
        configuration of each shape $\gamma$.  A natural family for
        this collection of pile configurations consists of what we
        call \emph{non-crossing pile configurations}; namely, for the
        composition $\gamma = (\gamma_{1}, \gamma_{2}, \ldots,
        \gamma_{k}) \models n$,
        \[
        \mathfrak{N} \cap
        \mathfrak{P}_{\gamma}(n) = \{\{ \{\gamma_{1} > \cdots > 1\},
        \{\gamma_{1} + \gamma_{2} > \cdots > \gamma_{1} + 1 \},
        \ldots, \{n > \cdots > n - \gamma_{k-1}\}\}\}
        \]
        so that there are exactly $2^{n-1}$ such pile configurations.
        One can also show that $\mathfrak{N}$ is the image
        $R(S_{n}(3\textrm{-}1\textrm{-}2))$ of all permutations
        avoiding the classical pattern $3\textrm{-}1\textrm{-}2$ under
        the Patience Sorting Algorithm.

    \section{Enumerating $S_{n}(3\textrm{-}\bar{1}\textrm{-}42)$}
    \label{sec:EnumeratingSn3142}

        In this section we use the results from Section
        \ref{sec:PileConfigurations} to both enumerate and
        characterize the permutations that avoid the generalized
        permutation pattern $2\textrm{-}31$ unless it is part of the
        generalized pattern $3\textrm{-}1\textrm{-}42$.  We call this
        restricted form of the generalized pattern $2\textrm{-}31$ a
        \emph{barred (generalized) permutation pattern} and denote it
        by $3\textrm{-}\bar{1}\textrm{-}42$.  (This notation is due to
        J. West, et al., and first appeared in the study of two-stack
        sortable permutations \cite{refDGG1998, refDGW1996, refWest1990}.)

         \begin{thm} ~\\[-20pt]
         \label{thm:EnumeratingSn3142}
            \begin{enumerate}
                \item The set of permutations
                $S_n(3\textrm{-}\bar{1}\textrm{-}42)$ that avoids the
                pattern $3\textrm{-}\bar{1}\textrm{-}42$ is exactly
                the set $RPW(R(\mathfrak{S}_{n}))$ of reverse patience
                words obtainable from the symmetric group
                $\mathfrak{S}_{n}$.\smallskip

                \item The size of
                $S_{n}(3\textrm{-}\bar{1}\textrm{-}42)$ is given by
                the $n^{\mathrm{th}}$ Bell number $B_{n}$.
            \end{enumerate}
        \end{thm}

        \begin{proof} ~\\[-20pt]
            \begin{enumerate}
                \item Let $\sigma \in
                S_{n}(3\textrm{-}\bar{1}\textrm{-}42)$.  Then, for $i
                = 1, 2, \ldots, n-1$, define $$\sigma_{m_{i}} = \min
                \{\sigma_{j} \ | \ i\leq j\leq n\}.$$ Since $\sigma$
                avoids $3\textrm{-}\bar{1}\textrm{-}42$, the
                subpermutation
                $\sigma_{i}\sigma_{i+1}\cdots\sigma_{m_{i}}$ must be a
                decreasing subsequence of $\sigma$.  (Otherwise
                $\sigma$ would necessarily contain a $2\textrm{-}31$
                pattern that is not part of a
                $3\textrm{-}1\textrm{-}42$ pattern.)  It follows that
                the left-to-right minima subsequences $s_{1}, s_{2},
                \ldots, s_{k}$ of $\sigma$ must be disjoint and
                satisfy Equation (\ref{eq:PileConfigurationCondition})
                so that the result follows by Lemmas
                \ref{lemma:PileConfiguration} and
                \ref{lem:ShadowDiagramPileCorrespondence}.

                \smallskip

                \item Recall that the Bell number $B_{n}$ enumerates
                the set partitions of the set $[n] = \{1, 2, \ldots,
                n\}$.  From Part (1), the elements of
                $S_{n}(3\textrm{-}\bar{1}\textrm{-}42)$ are in
                bijection with pile configurations.  Thus, since pile
                configurations are themselves set partitions by Lemma
                \ref{lemma:PileConfiguration}, we need only show that
                every set partition is also a pile configuration.  But
                this follows by ordering the components of a given set
                partition by their smallest element so that Equation
                (\ref{eq:PileConfigurationCondition}) is satisfied.
            \end{enumerate}
        ~\\[-35pt] %to get rid of extra space before end-of-proof box
        \end{proof}

        \begin{rem}
            We conclude by remarking that even though the set
            $S_{n}(3\textrm{-}\bar{1}\textrm{-}42)$ is enumerated by
            the very well known Bell numbers, it cannot be described
            in a simpler way using classical pattern avoidance.  This
            means that there does not exist a countable set of
            non-generalized (a.k.a. classical) permutation patterns
            $\tau_{1}, \tau_{2}, \ldots$ such that

            \begin{displaymath}
                    S_{n}(3\textrm{-}\bar{1}\textrm{-}42) =
                    S_{n}(\tau_{1}, \tau_{2}, \ldots) = \bigcap_{i
                    \geq 1} S_{n}(\tau_{i}).
            \end{displaymath}

            \noindent There are two very important reasons that this
            cannot happen:\medskip

            First of all, the Bell numbers satisfy $\log B_{n} = n
            (\log n - \log\log n + O(1))$ and so exhibit
            superexponential growth.  However, in light of the
            Stanley-Wilf ex-Conjecture (which was recently proven by
            Marcus and Tardos \cite{refMT2004}), the set of
            permutations $S_{n}(\tau)$ avoiding any classical pattern
            $\tau$ can only grow at most exponentially in $n$.

            On the other hand, the so-called \emph{avoidance class} of
            permutations
            \begin{displaymath}
                    Av(3\textrm{-}\bar{1}\textrm{-}42) = \bigcup_{n
                    \geq 0} S_{n}(3\textrm{-}\bar{1}\textrm{-}42)
            \end{displaymath}
            is not closed under taking order-isomorphic
            subpermutations, whereas it is easy to see that classes of
            permutations defined by classical pattern avoidance must
            be closed.  (See B\'{o}na \cite{refBona2004}, Chap.  5.)
            In particular, $3142 \in
            Av(3\textrm{-}\bar{1}\textrm{-}42)$ but $231 \notin
            Av(3\textrm{-}\bar{1}\textrm{-}42)$.

            \medskip

            At the same time, Theorem \ref{thm:EnumeratingSn3142}(2)
            implies that $3\textrm{-}\bar{1}\textrm{-}42$ belongs to
            the so-called \emph{Wilf Equivalence class} for the
            generalized pattern $1\textrm{-}23$.  That is, if $$\tau
            \in \{1\textrm{-}23, 3\textrm{-}21, 12\textrm{-}3,
            32\textrm{-}1, 1\textrm{-}32, 3\textrm{-}12,
            21\textrm{-}3, 23\textrm{-}1\}$$ then the size of the
            avoidance class $S_{n}(\tau)$ is also given by the
            $n^{\mathrm{th}}$ Bell number $B_{n}$.  In particular,
            Claesson \cite{refClaesson2001} showed that
            $|S_n(23\textrm{-}1)|=B_n$ via a direct bijection between
            permutations avoiding $23\textrm{-}1$ and set partitions.
            Furthermore, given $\sigma \in
            S_n(3\textrm{-}\bar{1}\textrm{-}42)$, each segment between
            consecutive right-to-left minima must be a single
            decreasing run (when from read left to right), so it is
            easy to see that
            $S_n(3\textrm{-}\bar{1}\textrm{-}42)=S_n(23\textrm{-}1)$.
            Thus, the barred pattern $3\textrm{-}\bar{1}\textrm{-}42$
            and the generalized pattern $23\textrm{-}1$ are not just
            in the same Wilf equivalence class but also have identical
            avoidance classes.

            Still, even though
            $S_n(3\textrm{-}\bar{1}\textrm{-}42)=S_n(23\textrm{-}1)$,
            it is more natural to use avoidance of
            $3\textrm{-}\bar{1}\textrm{-}42$ when studying Patience
            Sorting.  Fundamentally, this lets us look at
            $S_n(3\textrm{-}\bar{1}\textrm{-}42)$ as the set of
            equivalence classes in $\mathfrak{S}_{n}$ modulo
            $3\textrm{-}\bar{1}\textrm{-}42 \stackrel{PS}{\sim}
            3\textrm{-}\bar{1}\textrm{-}24$, where each equivalence
            class corresponds to a unique pile configuration.  The
            same equivalence relation is not easy to describe when
            starting with an occurrence of $23\textrm{-}1$.  (Note
            that $23\textrm{-}1 \sim 2\textrm{-}13$ or $23\textrm{-}1
            \sim 21\textrm{-}3$ is wrong since we would incorrectly
            get $2431 \sim 2314$ or $2431 \sim 2134$ instead of the
            correct $2431 \sim 2413$).

            \smallskip

            This suggests that there is even more information about
            pattern avoidance to be gotten from such a simple
            algorithm as Patience Sorting.

        \end{rem}

\end{document}